# Building Hyper Dirichlet Processes for Graphical Models

**Daniel Heinz**

*Carnegie Mellon University*
*5000 Forbes Avenue*
*Pittsburgh, PA 15213*
*e-mail:* `dheinz@stat.cmu.edu`

**Abstract:** Graphical models are used to describe the conditional independence relations in multivariate data. They have been used for a variety of problems, including log-linear models (Liu and Massam, 2006), network analysis (Holland and Leinhardt, 1981; Strauss and Ikeda, 1990; Wasserman and Pattison, 1996; Pattison and Wasserman, 1999; Robins et al., 1999), graphical Gaussian models (Roverato and Whittaker, 1998; Giudici and Green, 1999; Marrelec and Benali, 2006), and genetics (Dobra et al., 2004). A distribution that satisfies the conditional independence structure of a graph is *Markov*. A graphical model is a family of distributions that is restricted to be Markov with respect to a certain graph. In a Bayesian problem, one may specify a prior over the graphical model. Such a prior is called a *hyper Markov law* if the random marginals also satisfy the independence constraints. Previous work in this area includes (Dempster, 1972; Dawid and Lauritzen, 1993; Giudici and Green, 1999; Letac and Massam, 2007). We explore graphical models based on a non-parametric family of distributions, developed from Dirichlet processes.

## Contents

# 1. Introduction

Markov distributions are multivariate measures that satisfy a specified set of conditional independence relations. Often, an undirected graph is useful to represent this struture. A measure is Markov with respect to a graph if whenever two variables have no edge between them, they are conditionally independent given the other variables in the graph. Dawid and Lauritzen (1993) extended this notion to the parameter space. In Bayesian statistics, the measure of the data is random, and therefore has its own distribution called the prior. A prior law over Markov measures is *hyper Markov* if it gives probability one to Markov measures and the random marginal measures have the specified conditional independence structure. The hyper inverse Wishart distribution which serves as a prior for the multivariate Gaussian with known mean. The usual inverse Wishart is a specific case, which is hyper Markov for the saturated model.

Like all parametric models, the hyper inverse Wishart prior makes strong assumptions about the shape of the distribution. In many applications, such assumptions are undesirable. In contrast, non-parametric models make weak assumptions. Typical assumptions include continuity and the existence of some number of derivatives. For example, one may specify that the distribution is smooth, having derivatives of all orders. The current paper aims to apply the non-parametric approach to graphical models. We achieve this by following the framework laid by Dawid and Lauritzen. We begin with the Dirichlet Process, a commonly used non-parametric prior law. We then describe how to build this family into a non-parametric hyper Markov prior.

As in Dawid and Lauritzen (1993), we restrict our attention to decomposable graphs. The benefit of this is that a decomposable graph can be easily built up from smaller components called cliques which intersect to form the entire graph. Dawid and Lauritzen begin by considering a base distribution for each clique. The only requirement is that these distributions agree where the cliques intersect. They weave together base distributions by taking the base measure of one clique as its marginal, and conditioning the second clique on the intersection. They repeat this process for every clique. The third clique is added by conditioning its base measure on the intersection with the previous two cliques. This process is repeated until all the cliques have been combined. The end result is a Markov distribution whose marginal over each clique is the clique's base. For a prior on Markov distributions, Dawid and Lauritzen construct a hyper Markov law in the same way.

As an example of the Dawid and Lauritzen (1993) construction, consider the problem of estimating the covariance matrix of a multivariate Gaussian. If we believe that the data exhibit some conditional independence structure, this implies certain constraints on the covariance matrix. (Speed and Kiiveri, 1986) showed that the sufficient statistics are the component covariance matrices belonging to each clique. The inverse Wishart is the usual prior for the saturated model which has no constraints on the covariance matrix. In a non-saturated model, the sub-matrix of each clique is unconstrained, except that the sub-matrices must agree where their indices intersect. For this reason, the inverse

Wishart is the natural choice as the base measure for each clique. The sub-matrix for the first clique has an inverse Wishart prior. If the graph is connected and the cliques have a perfect ordering (see Section 2.1), then the first and second sub-matrices have some elements in common. Thus, the sub-matrix for the second clique is the inverse Wishart, conditional on knowing some of the elements. By repeating the conditioning for each clique, the *hyper inverse Wishart* is defined.

In the current paper, we apply this framework to non-parametric priors. Instead of the inverse Wishart, the Dirichlet process prior is the base measure for each clique. Following the analogy, we build the marginals into a hyper Markov prior, which we refer to as the hyper Dirichlet process. The Dirichlet process is a special case of tail-free processes (Ferguson, 1973). Dirichlet processes have been used for non-parametric priors in many areas, including block models (Bush and MacEachern, 1996), survival analysis (Susarla and Ryzin, 1976; Ghosh and Ramamoorthi, 1995; Kim and Lee, 2001), and non-stationary point processes (Pievatolo and Rotondi, 2000). These are all areas that could potentially use a hyper Dirichlet process in multidimensional problems. In section 2, we explain notation and formalize some of the ideas presented so far. In section 3, we describe the Dirichlet and some previous results. In section 4 we weave Dirichlet processes on the cliques to build the hyper Dirichlet process, and show that it is a hyper Markov prior. Finally, we explore applications for this framework in section 5.

## 2. Notation and Setting

Throughout this paper we consider a graph, $\mathcal{G}$, with vertex set $\mathbf{V}$ and edge set $\mathbf{E}$. By convention, we assume that $(\gamma, \gamma) \in \mathbf{E}$ for all $\gamma$. We call such edges *loops*. There is no practical difference if loops are excluded from $\mathbf{E}$, though some minor changes are required for certain definitions. If $\mathbf{A} \subseteq \mathbf{V}$, then $\mathcal{G}_\mathbf{A}$ is the subgraph of $\mathcal{G}$ over $\mathbf{A}$. The subgraph $\mathcal{G}_\mathbf{A}$ has vertex set $\mathbf{A}$, and edge set $\mathbf{E}_\mathbf{A} = (\mathbf{A} \times \mathbf{A}) \cap \mathbf{E}$. We say that $\mathbf{A}$ induces the subgraph $\mathcal{G}_\mathbf{A}$. If $\mathbf{E}_\mathbf{A} = \mathbf{A} \times \mathbf{A}$, then $\mathcal{G}_\mathbf{A}$ is complete. A *clique* is a set $\mathbf{A}$ such that $\mathcal{G}_\mathbf{A}$ is complete and for any superset $\mathbf{B} \supset \mathbf{A}$, $\mathcal{G}_\mathbf{B}$ is not complete. For example, if $\mathcal{G}$ itself is complete, then there is one clique, viz. $\mathbf{V}$.

A $k$-path is a sequence $(\gamma_0, \gamma_1, \ldots, \gamma_k)$, such that $(\gamma_i, \gamma_{i+1}) \in \mathbf{E}$, for $0 \leq i < k$. If $\mathbf{A}$ and $\mathbf{B}$ are subsets of $\mathbf{V}$, then a path between them is a path between any $a \in \mathbf{A}$ and any $b \in \mathbf{B}$. A graph is *connected* if there exists a path between every pair of subsets. A third subset $\mathbf{C} \in \mathbf{V}$ is said to *separate* $\mathbf{A}$ and $\mathbf{B}$ if every path between them contains an element of $\mathbf{C}$. A $k$-cycle is a path such that $k \geq 3$, $\gamma_0 = \gamma_k$ and the other elements are distinct. A graph is *decomposable* if there are no cycles longer than length 3. A decomposable graph admits a *perfect ordering* of its cliques.

**Definition 1** PERFECT ORDERING. *Suppose a graph $\mathcal{G}$ has $n$ cliques. Let the cliques have an arbitrary ordering $\mathbf{C}_1, \ldots, \mathbf{C}_n$. Define $\mathbf{H}_k = \cup_{i=1}^{k} \mathbf{C}_i$. For $k \geq 2$ define $\mathbf{S}_k = \mathbf{C}_k \cap H_{k-1}$ and $\mathbf{R}_k = \mathbf{C}_k \setminus H_{k-1}$. The ordering of the cliques is a perfect ordering if for each $2 \leq k \leq n$, there exists $j_k < k$ such that $\mathbf{S}_k \subset \mathbf{C}_{j_k}$.*

The sets $\mathbf{H}_k$ are called the histories. The separators, $\mathbf{S}_k$, separate $\mathbf{C}_k$ from the previous history. The sets $\mathbf{R}_k$ are called the residuals, which represent the new nodes being added to the history. In a perfect ordering, each new clique is separated from the current set of nodes by a single one of the earlier cliques.

For every, $\gamma \in \mathbf{V}$, $X_\gamma$ is a random variable taking values in the space $(\mathcal{X}_\gamma, \mathcal{F}_\gamma)$. In this sense, we consider $\mathbf{V}$ an index set of components of some random variable $X = (X_\gamma : \gamma \in \mathbf{V})$. We denote the range and $\sigma$-field of $X$ by $(\mathcal{X}, \mathcal{F}) = (\times_{\gamma \in \mathbf{V}} \mathcal{X}_\gamma, \times_{\gamma \in \mathbf{V}} \mathcal{F}_\gamma)$. Furthermore, we extend these definitions to subsets, $\mathbf{A} \subseteq \mathbf{V}$.

$$\begin{aligned} X_\mathbf{A} &= (X_\mathbf{A} : \gamma \in \mathbf{A}) \\ \mathcal{X}_\mathbf{A} &= \times_{\gamma \in \mathbf{A}} \mathcal{X}_\gamma \\ \mathcal{F}_\mathbf{A} &= \times_{\gamma \in \mathbf{A}} \mathcal{F}_\gamma \end{aligned}$$

Let $\alpha$ be a measure over some $\mathcal{X}_\mathbf{A}$, then $\overline{\alpha} = \alpha/\alpha(\mathcal{X}_\mathbf{A})$. In other words, $\overline{\alpha}$ is the probability measure proportional to $\alpha$. If $\mathbf{B} \subseteq \mathbf{A}$, then $\alpha_\mathbf{B}$ is the marginal of $\alpha$ over $\mathcal{X}_\mathbf{B}$. Thus, $\alpha_\mathbf{B}(U) = \alpha(U \times \mathcal{X}_{\mathbf{A} \setminus \mathbf{B}})$, $\forall\, U \in \mathcal{F}_\mathbf{B}$. For convenience, if $\gamma \in \mathbf{A}$, we may write $\theta_\gamma$ for $\theta_{\{\gamma\}}$. If $\alpha$ and $\beta$ are both measures on some space $(\mathcal{X}, \mathcal{F})$, then we define their sum, $\alpha + \beta$, by

$$[\alpha + \beta](U) = \alpha(U) + \beta(U), \quad \forall\, U \in \mathcal{F}.$$

If $x \in \mathcal{X}$, then the delta measure $\delta_x$ is a point mass concentrated at $x$:

$$\delta_x(U) = \begin{cases} 1, & x \in U \\ 0, & x \notin U \end{cases}, \quad \forall U \in \mathcal{F}.$$

### 2.1. Graph Selection

For the remainder of the paper, we consider undirected graphs, which implies that $(i,j) \in \mathbf{E}$ if and only if $(j,i) \in \mathbf{E}$. We also assume that the graph is connected and decomposable. An undirected graph depicts the conditional independence structure for some variable $X$ in the following sense:

$$X_\mathbf{A} \perp\!\!\!\perp X_\mathbf{B} \mid X_\mathbf{C} \text{ whenever } \mathbf{C} \text{ separates } \mathbf{A} \text{ and } \mathbf{B}. \tag{1}$$

**Definition 2** MARKOV PROBABILITY MEASURE. *If $\theta$ is a probability measure on $(\mathcal{X}, \mathcal{F})$, we say it is Markov on a decomposable graph, $\mathcal{G}$, if $X \sim \theta$ satisfies the conditional independences in $\mathcal{G}$.*

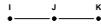

FIG 1. *A graph depicting conditional independence of I and K given J.*

**Example 1** *Let $\mathcal{G}$ be the graph depicted in Figure 2.1. A measure $\theta$ is Markov on $\mathcal{G}$, if and only if $X_I \perp\!\!\!\perp X_K | X_J$ whenever $X \sim \theta$.*

Implicit in the definition is the fact that it is only sensible to refer to a measure as Markov in relation to a specific graph. For example, if the measure $\theta$ is not Markov on $\mathcal{G}$ in Example 1, it is still Markov on the saturated graph with $\mathbf{V} = \{I, J, K\}$. All measures over $\mathcal{X}_\mathbf{V}$ are trivially Markov on the saturated graph since there are no constraints on conditional independence. Furthermore, if $\mu$ is a measure such that each $X_\gamma$ is independent, then it is Markov on any graph (with the appropriate vertex set.) We denote the set of all distributions that are Markov on $\mathcal{G}$ by $\mathscr{M}(\mathcal{G})$.

It will be useful to keep Figure 1 in mind throughout this paper. While the graph technically has only three variables, it is representative of any connected graph of two cliques. Instead of one variable, imagine $I$, $J$, and $K$ to contain multiple variables, with $J$ being the variables that belong to both cliques. $I$ is the set of variables in one clique but not the other, and $K$ vice versa.

Let $X \sim P \in \mathscr{F}$ be a random variable whose distribution is modeled by some family of probability distributions. In some applications, the focus is not on determining $P$, but on discovering the independence structure of $X$. A graph of this structure, $\mathcal{G}$, denotes the belief that $P$ is Markov with respect to $\mathcal{G}$. Thus, it restricts the model to a sub-family, $\mathscr{F}_\mathcal{G} = \mathscr{F} \cap \mathscr{M}(\mathcal{G})$. Graph selection is the problem of determining the smallest $\mathscr{F}_\mathcal{G}$ that contains $P$. The most prevalent examples are graphical Gaussian models. Graph selection for Gaussian models is often called covariance selection. In this setting, the relevant family is the set of $p$-variate Gaussian distributions. Denote this family $\mathscr{N} = \{N_p(\mu, \Sigma) : \mu \in \mathcal{R}^p, \Sigma \in M_p^+\}$, where $M_p^+$ is the cone of real-valued, symmetric $p \times p$ matrices that are positive definite. Specifying a graph, $\mathcal{G}$, translates to putting constraints on $\Sigma$. For example, if $(x_1, x_2, x_3)$ is such that $x_1 \perp\!\!\!\perp x_3 | x_2$, then $\sigma_{13}$ is no longer a free parameter, but a function of $\sigma_{11}, \sigma_{22}, \sigma_{33}, \sigma_{12}$, and $\sigma_{23}$. In general, denote the sub-family of Gaussian distributions Markov on $\mathcal{G}$ by $\mathscr{N}_\mathcal{G}$. Let $P_\mathcal{G}$ be the set of positive definite matrices such that $K_{ij} = 0$ for all $(i, j) \notin \mathbf{E}$; let $Q_\mathcal{G}$ be the image of $P_\mathcal{G}$ under matrix inversion. Speed and Kiiveri (1986) showed that if $N_p(\mu, \Sigma_\mathcal{G})$ is Markov with respect to $\mathcal{G}$, then $\Sigma_\mathcal{G} \in Q_\mathcal{G}$. Thus $\mathscr{N}_\mathcal{G} = \{N_p(\mu, \Sigma) : \mu \in \mathcal{R}^p, \Sigma \in Q_\mathcal{G}\}$. The goal of covariance selection is to find the smallest $Q_\mathcal{G}$ containing $\Sigma$, the population covariance matrix.

Much progress has been made with graph selection for parametric models. Dawid and Lauritzen (1993) proved many results for decomposable graphical models, including multinomial and multivariate Gaussian problems. For example, they present the distribution of the restricted maximum likelihood estimate of $\Sigma$ for the $\mathscr{N}_\mathcal{G}$ model with $\mu$ known. This distribution is called the hyper Wishart distribution since it reduces to the Wishart when $\mathcal{G}$ is complete. Letac and Massam (2007) have extended the hyper Wishart to a richer family of distributions on $Q_\mathcal{G}$ and $P_\mathcal{G}$. Giudici and Green (1999) implemented a reversible jump Markov chain Monte Carlo algorithm for determining $\mathcal{G}$.

The family of hyper inverse Wishart distributions is the subset of Markov distributions such that each clique marginal is inverse Wishart. Carvalho et al. (2007) provide an algorithm for generating random variables from this family. For decomposable models, the presence of a perfectordering simplifies the process. For two cliques, the algorithm begins by generating an inverse Wishart variable on one clique. If the cliques overlap, this determines some of the parameters for the other clique. Therefore, one needs to generate a conditional Wishart variable given those entries. For multiple cliques, one simply repeats this process. With a perfect ordering, the process is simplified because each new clique is conditioned on only one previous clique. Conditioning on multiple cliques can lead to moderate complications in the conditional distribution. Hence, decomposable models are computationally convenient.

## 3. The Dirichlet Process

The Dirichlet process (Ferguson, 1973) is a special case of tail-free distributions. It is a prior, meaning that it provides a distribution over the space of probability distributions on $(\mathcal{X}, \mathcal{F})$. In this paper, we use the term *law* to refer to distributions over probability measures. However, this terminology is merely a convenience; the words "law" and "distribution" are typically interchangeable. The Dirichlet process is an example of a non-parametric law, which means that it cannot be specified by a finite-dimensional parameter. In this section, the Dirichlet process is introduced and some of its useful properties are given. This leads into the next section, in which we show how to compose a hyper Dirichlet processes from multiple Dirichlet processes.

**Definition 3** DIRICHLET PROCESS. *Let $\mathbf{A}$ be any subset of $\mathbf{V}$. Let $\alpha$ be a measure over $(\mathcal{X}_\mathbf{A}, \mathcal{F}_\mathbf{A})$, and let $\theta$ be a random probability measure over the same space. We say that the distribution of $\theta$ is a Dirichlet process with base measure $\alpha$, and write $\theta \sim DP_\alpha$, if*

$$(\mathbb{P}(A_1), \mathbb{P}(A_2), \ldots, \mathbb{P}(A_k)) \sim \mathrm{Dir}(\alpha(A_1), \alpha(A_2), \ldots, \alpha(A_k)), \qquad (2)$$

*whenever $(A_i)_{i=1}^k$ is a partition of $\mathbf{A}$.*

This definition leads to some useful properties, described in the following theorem.

**Theorem 4** POSTERIOR DIRICHLET PROCESS. *Let $\theta \sim DP_\alpha$ and, given $\theta$, let $X_1, \ldots, X_n$ be an iid sample from $\theta$.*

*(i) $X_i \sim \overline{\alpha} \ \forall \ i$.*
*(ii) $\theta | (X_1, \ldots, X_n) \sim DP_{\alpha'}$, where $\alpha' = \alpha + \sum_{i=1}^n \delta_{X_i}$.*

See Theorem 1.9.4 of Schervish (1995), p. 54. □

The first property states that if the random measure is integrated out, the marginal distribution of the data is $\overline{\alpha}$. This property ensures that a Markov base measure implies that the Dirichlet process, integrated over all possible $\theta$, is a Markov distribution. This does not guarantee that the process is a hyper Markov law. That requires the stronger condition that $\theta \sim DP_\alpha$ is a Markov distribution with probability one. The second property states that if a Dirichlet process is used as a prior measure, then the posterior measure is also a Dirichlet process, with an easily updated base measure. This fact helps determine which properties of a prior will persist in the posterior.

If the prior law of $\theta$ is a Dirichlet process, then the various marginal distributions of $\theta$ will also have a Dirichlet process law. This is expressed in the following theorem.

**Theorem 5** MARGINAL OF A DIRICHLET PROCESS. *Let $\theta \sim DP_\alpha$ be a random probability measure on $(\mathcal{X}_\mathbf{A}, \mathcal{F}_\mathbf{A})$. For $\mathbf{B} \subseteq \mathbf{A}$, the marginal of $\theta$ over $\mathbf{B}$ is $\theta_\mathbf{B} \sim DP_{\alpha_\mathbf{B}}$.*

PROOF: Define $\mathbf{A}' = \mathbf{A} \setminus \mathbf{B}$. Let $B_1, B_2, \ldots, B_k$ be a measurable partition of $\mathbf{B}$.

$$\begin{aligned}
(\mathbb{P}_{\theta_\mathbf{B}}(B_1), \ldots, \mathbb{P}_{\theta_\mathbf{B}}(B_k)) &= (\mathbb{P}_\theta(B_1 \times \mathcal{X}_{\mathbf{A}'}), \ldots, \mathbb{P}_\theta(B_k \times \mathcal{X}_{\mathbf{A}'})) \\
&\sim \text{Dir}(\alpha(B_1 \times \mathcal{X}_{\mathbf{A}'}), \ldots, \alpha(B_k \times \mathcal{X}_{\mathbf{A}'})) \\
&= \text{Dir}(\alpha_\mathbf{B}(B_1), \ldots, \alpha_\mathbf{B}(B_k)). \qquad \square
\end{aligned}$$

We proceed by showing how the Dirichlet process can be used as a nonparametric prior (see Ferguson (1973) for details.) Let $F$ be an unknown cumulative probability distribution that we wish to estimate. For simplicity, we consider a one-dimensional random variable. Let $\pi = DP_\alpha$ be the prior law. Let the loss function be a squared error loss. Then the Bayes' risk is

$$R_\pi(F, \widetilde{F}) = \int E(F(t) - \widetilde{F}(t))^2 dt. \tag{3}$$

The risk is minimized by setting $\widetilde{F}(t)$ to $EF(t)$, where the expectation is relative to the posterior distribution. If we observe data $X_1, X_2, \ldots X_n$, then the posterior is $DP_{\alpha'}$, where $\alpha' = \alpha + \sum_i \delta_{X_i}$ (see Theorem 7.) The posterior distribution of $P((-\infty, t])$ is $\text{Beta}(\alpha'((-\infty, t]), \alpha'((t, \infty)))$. Therefore,

$$\begin{aligned}
EF(t) &= \frac{\alpha'((-\infty, t])}{\alpha'((-\infty, t]) + \alpha'((t, \infty))} \tag{4} \\
&= \frac{\alpha((-\infty, t]) + \sum_{i=1}^n 1_{(X_i \leq t)}}{\alpha(\mathcal{X}) + n}. \tag{5}
\end{aligned}$$

The Bayes estimate can be written as a weighted sum of two estimates

$$\begin{aligned}
\widetilde{F}(t) &= EF(t) \tag{6} \\
&= (1-w)\overline{\alpha}((-\infty, t]) + w\widehat{F}(t), \tag{7}
\end{aligned}$$

where $\overline{\alpha}((-\infty, t])$ is the prior estimate, $\widehat{F}(t)$ is the empirical cdf, and $w = n(\alpha(\mathcal{X}) + n)^{-1}$ is the weight of the data. This convex combination of a prior estimate and frequentist estimate is common in Bayesian analysis. This shows the role of the base measure in the Dirichlet process. $\overline{\alpha}$ is the prior guess about the shape of the unknown distribution. $\alpha(\mathcal{X})$ is mathematically equivalent to the prior sample size.

### 3.1. Dirichlet Process as a Stick-Breaking Prior

A stick-breaking process is an almost surely discrete random probability measure, $\theta$, that can be expressed as

$$\theta(\cdot) = \sum_{k=1}^{N} w_k \delta_{Z_k}(\cdot), \qquad (8)$$

where the $Z_k$ are independently distributed atoms from some distribution $G$, and $\sum_{k=1}^{N} w_k = 1$ almost surely. The number of atoms, $N$, may be finite or infinite. The weights are determined by successively breaking random pieces of a unit-length stick. Thus, $w_1 = p_1$, $w_2 = (1-p_1)p_2$, and $w_k = p_k \prod_{i=1}^{k-1}(1-p_i)$. For finite $N$, $w_N$ is defined by $1 - \sum_{i=1}^{N-1} w_i$, or equivalently by $\prod_{i=1}^{N-1}(1-p_i)$. Traditionally, stick-breaking measures are defined such that $p_k$ is defined as a Beta$(a_k, b_k)$ random variable for $1 \leq k < N$. Thus, a stick-breaking measure is specified by a probability distribution $P$, and a countable sequence of Beta parameters $(a_k, b_k)_{k=1}^{N-1}$. Sethuraman (1994) showed that a Dirichlet Process is a stick-breaking measure with $Z_k \sim \overline{\alpha}$, and $(a_k, b_k) = (0, \alpha(\mathcal{X}))$ for all $k \in \mathbb{N}$. This relationship leads to an alternative definition of the Dirichlet process.

**Definition 6** DIRICHLET PROCESS (alternate definition). *Let $\mathbf{A}$ be any subset of $\mathbf{V}$. Let $G$ be a probability measure on $(\mathcal{X}_\mathbf{A}, \mathcal{F}_\mathbf{A})$, and let $\theta$ be a random probability measure over the same space. For $\nu > 0$, we say that the distribution of $\theta$ is a Dirichlet process with base distribution (or measure) $G$ and precision $\nu$, and write $\theta \sim DP(\nu G)$, if*

$$(\mathbb{P}(A_1), \mathbb{P}(A_2), \ldots, \mathbb{P}(A_k)) \sim \text{Dir}(\nu G(A_1), \nu G(A_2), \ldots, \nu G(A_k)), \qquad (9)$$

*whenever $(A_i)_{i=1}^{k}$ is a partition of $\mathbf{A}$.*

Note that this distribution is equivalent to Definition 3 by letting $\alpha = \nu G$. For example, $\nu$ is equivalent to the prior sample size, and $G$ is equivalent to the prior mean. In this definition, $\nu$ and $G$ are easily translated as the parameters of a stick-breaking measure. That is, the random atoms are iid $G$, and $p_k \sim$ Beta$(0, \nu)$ for all $k \in \mathbb{N}$. Because the stick-breaking representation is useful for many of the theorems we prove, Definition 6 will be the definition of choice for much of the current paper.

The previous theorems regarding Dirichlet processes can be expressed using $\nu G$ notation. For example, we rewrite Theorem 4 regarding the posterior Dirichlet process.

**Theorem 7** POSTERIOR DIRICHLET PROCESS (alternate). *Let $\theta \sim DP(\nu G)$ and, given $\theta$, let $X_1, \ldots, X_n$ be an iid sample from $\theta$.*

(i) $X_i \sim G \ \forall \ i$.
(ii) $\theta|(X_1, \ldots, X_n) \sim DP_{\nu'G'}$, where $\nu' = \nu + n$ and $G' = (\nu + n)^{-1}(\nu P + \sum_{i=1}^{n} \delta_{X_i})$.

In the following section we introduce the hyper Dirichlet process and show that it is an example of a stick-breaking measure. We then use Equation 8 to prove some of its properties. While we focus on the hyper Dirichlet Process for simplicity and concreteness, many of the results apply to other stick-breaking processes as well.

## 4. The Hyper Dirichlet Process

Consider a multivariate variable $X$ with distribution $\theta$. Suppose that we know little about $\theta$, other than it is Markov on some decomposable graph, $\mathcal{G}$. In this case we may wish to specify a non-parametric prior for $\theta$. For example, we focus on the Dirichlet process. There are two main difficulties with this approach. The first is the elicitation of a proper base measure. The second is ensuring that the Dirichlet process gives probability one to $\mathcal{M}(\mathcal{G})$. Both concerns are addressed by using a framework that we dub the *hyper Dirichlet process*.

To define a hyper Dirichlet process, we begin by eliciting a base measure for each clique in $\mathcal{G}$. Hopefully, this is simpler than eliciting a base measure for the entire graph at once. These base measures are combined to form a base measure over the entire graph. We define these combinations in a way which will ensure that the support of the process lies within the set of Markov distributions on $\mathcal{G}$. In the remainder of this section, we provide details to this method and show that it satisfies the Markov property.

### 4.1. Markov Combinations of Probability Measures

Dawid and Lauritzen (1993) show that if two subsets of $\mathbf{V}$ are each endowed with a marginal probability measure, then there is a logical choice for their joint distribution, provided the marginals satisfy a consistency condition.

**Definition 8** CONSISTENCY (OF PROBABILITY MEASURES). *Suppose* $\mathbf{A}, \mathbf{B} \subseteq \mathbf{V}$. *Let $\mu$ and $\lambda$ be probability measures on $(\mathcal{X}_\mathbf{A}, \mathcal{F}_\mathbf{A})$ and $(\mathcal{X}_\mathbf{B}, \mathcal{F}_\mathbf{B})$, respectively. We say that $\mu$ and $\lambda$ are consistent if they induce the same marginal over $\mathcal{X}_{\mathbf{A} \cap \mathbf{B}}$.*

Note that $\mu$ and $\lambda$ are consistent only if

$$\mu(\mathcal{X}_{\mathbf{A} \setminus \mathbf{B}} \times U) = \lambda(\mathcal{X}_{\mathbf{B} \setminus \mathbf{A}} \times U) \ \ \forall \ U \in \mathcal{F}_{\mathbf{A} \cap \mathbf{B}}.$$

**Theorem 9** *Suppose $\mu$ on $(\mathcal{X}_\mathbf{A}, \mathcal{F}_\mathbf{A})$ and $\lambda$ on $(\mathcal{X}_\mathbf{B}, \mathcal{F}_\mathbf{B})$ are consistent probability measures, with $\mathbf{A}, \mathbf{B} \subseteq \mathbf{V}$. There exists an almost-everywhere unique distribution, $\alpha$, such that:*

(i) $\alpha_\mathbf{A} = \mu$,
(ii) $\alpha_\mathbf{B} = \lambda$.

*(iii)* $\alpha \in \mathscr{M}(\mathcal{G}_{\mathbf{A} \cup \mathbf{B}})$,

PROOF: Construct $\alpha$ such that its marginal over $\mathcal{X}_{\mathbf{A}}$ is $\mu$, so that condition (i) is satisfied. Specify its conditional distributions over $\mathcal{X}_{\mathbf{B}}$ given $X_{\mathbf{A}}$ to be the same as the conditional distributions of $\lambda$ given $X_{\mathbf{A} \cap \mathbf{B}}$. This ensures that (iii) holds as well. Let $\mathbf{C} = \mathbf{A} \cap \mathbf{B}$ and $\mathbf{B}' = \mathbf{B} \setminus \mathbf{A}$. Then for any $U \in \mathbf{B}'$ and $V \in \mathbf{C}$,

$$\begin{aligned}
\mathbb{P}_{\alpha_{\mathbf{B}}}(U \times V) &= \mathbb{P}_{\alpha_{\mathbf{B}'|\mathbf{C}}}(U|V)\mathbb{P}_{\alpha_{\mathbf{C}}}(V) \\
&= \mathbb{P}_{\lambda_{\mathbf{B}'|\mathbf{C}}}(U|V)\mathbb{P}_{\mu_{\mathbf{C}}}(V) \\
&= \mathbb{P}_{\lambda_{\mathbf{B}'|\mathbf{C}}}(U|V)\mathbb{P}_{\lambda_{\mathbf{C}}}(V) \\
&= \mathbb{P}_{\lambda}(U \times V).
\end{aligned}$$

The second equation follows from the construction of $\alpha$. The third equation is ensured since $\mu$ and $\lambda$ are consistent. Hence, condition (ii) is also satisfied. Furthermore, the conditional distributions are unique, except over some subset of $\mathcal{X}_{\mathbf{C}}$ with zero measure under $\lambda$, and hence also under $\mu$ by consistency. Therefore, this construction gives (a version of) the unique distribution satisfying the conditions.

**Definition 10** MARKOV COMBINATION (OF PROBABILITY MEASURES). *Let $\mu$ and $\lambda$ be as in Theorem 9. We call the unique distribution satisfying (i)-(iii) the Markov Combination of $\mu$ and $\lambda$, and denote it by $\mu \star \lambda$.*

Now suppose $\mathcal{G}$ has a perfect ordering of cliques $(\mathbf{C}_1, \mathbf{C}_2, \ldots, \mathbf{C}_k)$, and that each clique $\mathbf{C}_i$ is imbued with a marginal probability distribution $P_i$. Further suppose that $P_i$ and $P_j$ are consistent for all $i, j$. Each clique is consistent with the previous history regarding the separator, since the separator is contained by a single previous clique. Using the idea of a Markov combination iteratively, we stitch together a distribution that is Markov on $\mathcal{G}$ and has the given marginals. Define $G_1 = P_1$, and $G_i = G_{i-1} \star P_i$ for $i \geq 2$. Dawid and Lauritzen show that $G = G_k$ is the unique Markov distribution satisfying $G_{\mathbf{C}_k} = P_k$. We call $G$ the Markov combination of $P_1, \ldots, P_k$. In general, we may write $\star(Q_1, \ldots, Q_k)$ to indicate a Markov combination with the understanding that the cliques are perfectly ordered and $Q_1, \ldots, Q_k$ are pairwise consistent.

### 4.2. Markov Combinations of Finite Measures

Using Markov combinations, we are able to take probability distributions and build a distribution over the entire graph. The base measure of a Dirichlet process, however, is not necessarily a probability distribution. Therefore, we proceed by extending Markov combinations to finite measures in general. For probability measures, we required the conditionals over $\mathcal{X}_{\mathbf{A} \cap \mathbf{B}}$ to be the same. We simply extend this definition to any finite measure.

**Definition 11** CONSISTENCY OF FINITE MEASURES. *Let $\mu$ be a finite measure over $(\mathcal{X}_{\mathbf{A}}, \mathcal{F}_{\mathbf{A}})$ and $\lambda$ be a finite measure over $(\mathcal{X}_{\mathbf{B}}, \mathcal{F}_{\mathbf{A}})$. We say that $\mu$ and $\lambda$* are consistent if they induce the same marginal measure over $\mathbf{A} \cap \mathbf{B}$. That is, $\mu$ and $\lambda$ are consistent if

$$\mu(\mathcal{X}_{\mathbf{A} \setminus \mathbf{B}} \times U) = \lambda(\mathcal{X}_{\mathbf{A} \setminus \mathbf{B}} \times U) \quad \forall \ U \in \mathcal{F}_{\mathbf{A} \cap \mathbf{B}}. \tag{10}$$

Recall that $\overline{\mu}$ is the probability measure proportional to $\mu$. Equation 10 holds if the following two conditions are satisfied:

1. $\overline{\mu}$ and $\overline{\lambda}$ are consistent.
2. $\mu(X_{\mathbf{A}}) = \lambda(X_{\mathbf{B}})$.

Consider these two conditions in the context of base measures for Dirichlet processes. $\overline{\mu}$ is the prior guess about the probability distribution of $X_{\mathbf{A}}$, and $\overline{\lambda}$ is the prior guess for $X_{\mathbf{B}}$. The first condition therefore states that the priors must agree about the distribution of $X_{\mathbf{A} \cap \mathbf{B}}$. It is reasonable to require that our prior is coherent in this way. The second condition states that the prior sample sizes for both sets of variables must be equal. This restraint is perhaps less desirable. It would be perfectly logical to be more certain about certain dimensions than others. Unfortunately, any measure on $\mathcal{X}_{\mathbf{A} \cup \mathbf{B}}$ must satisfy

$$\alpha_{\mathbf{A}}(\mathcal{X}_{\mathbf{A}}) = \int_{\mathcal{X}_{\mathbf{A}}} d\alpha_{\mathbf{A}} = \int_{\mathcal{X}_{\mathbf{A}}} \int_{\mathcal{X}_{\mathbf{B} \setminus \mathbf{C}}} d\alpha = \alpha(\mathcal{X}_{\mathbf{A} \cup \mathbf{B}}) = \int_{\mathcal{X}_{\mathbf{B}}} \int_{\mathcal{X}_{\mathbf{B} \setminus \mathbf{A}}} d\alpha = \int_{\mathcal{X}_{\mathbf{B}}} d\alpha_{\mathbf{B}} = \alpha_{\mathbf{B}}(\mathcal{X}_{\mathbf{B}}). \tag{11}$$

If $\mu(\mathcal{X}_{\mathbf{A}}) \neq \lambda(\mathcal{X}_{\mathbf{B}})$ there is no measure $\alpha$ on $\mathcal{X}_{\mathbf{A} \cup \mathbf{B}}$ satisfying $\alpha_{\mathbf{A}} = \mu$ and $\alpha_{\mathbf{B}} = \lambda$. In some situations, this problem is not too severe. Using the alternative definition, we express $\mu = \nu_1 G_1$ and $\lambda = \nu_2 G_2$. The consistency conditions translate to $G_1 = G_2$ and $\nu_1 = \nu_2$. If only the second condition fails, then it is still possible to find $G = G_1 \star G_2$. Employing the stick-breaking condition, we can generate random atoms from $G$. The problem lies in assigning weights to each atom. However, we can take solace in a mitigating factor. For density estimation, the value of the prior precision ($\nu$) is typically small compared to the sample size. Hence, we may be reasonable to simply scale the precisions so that $\nu_1$ and $\nu_2$ are equal. For these applications, it is only important that $G_1$ and $G_2$ are consistent. If so, the base measures $\mu$ and $\lambda$ only need to be *proportional* to each other over $\mathbf{A} \cap \mathbf{B}$.

There may be other situations in which scale *is* important. Unfortunately, as Equation 11 shows, we cannot find a suitable base measure for the prior that satisfies both $\mu$ and $\lambda$. Without a suitable prior, there can be no suitable posterior. If the goal is to estimate a distribution and there is genuine concern about the precision of the prior estimate, then both conditions must be satisfied. That is, $\mu$ and $\lambda$ must be consistent. Equivalently, $G_1$ and $G_2$ must be consistent and $\nu_1$ must be equal to $\nu_1$. Cases in which one or both conditions fail are explored more fully in Appendix A.

Subsequently, we assume that both consistency conditions are satisfied. This leads to a natural extension of the previous work. We have equated consistency of base measures with consistency of probability measures. Thus, we generalize Markov combinations to include consistent finite measureszed by scaling them to probability measures, finding the Markov combination, and rescaling the measures.

**Definition 12** MARKOV COMBINATION OF FINITE MEASURES. *Let $\mu$ be a finite measure on $(\mathcal{X}_\mathbf{A}, \mathcal{F}_\mathbf{A})$. Let $\lambda$ be a finite measure on $(\mathcal{X}_\mathbf{B}, \mathcal{F}_\mathbf{B})$ that is consistent with $\mu$. The Markov combination of $\mu$ and $\lambda$ is denoted $\mu \star \lambda$, where*

$$\mu \star \lambda = \mu(\mathcal{X}_\mathbf{A}) \cdot [\overline{\mu} \star \overline{\lambda}] = \lambda(\mathcal{X}_\mathbf{B}) \cdot [\overline{\mu} \star \overline{\lambda}], \tag{12}$$

*where $[\overline{\mu} \star \overline{\lambda}]$ is the almost-everywhere unique probability distribution satisfying Theorem 9.*

This definition is a generalization of Definition 10 for probability measures. Note that the Markov combination defined in this way is unique almost everywhere, since $[\overline{\mu} \star \overline{\lambda}]$ is unique almost everywhere.

It is easy to show that the $\overline{\phantom{x}}$ and $\star$ operations commute (with respect to composition).

**Theorem 13** *If $\mu$ and $\lambda$ are consistent measures, then $\overline{\mu \star \lambda} = \overline{\mu} \star \overline{\lambda}$.*

PROOF:

$$\overline{\mu \star \lambda} = \frac{[\mu \star \lambda]}{[\mu \star \lambda](\mathcal{X}_{\mathbf{A} \cup \mathbf{B}})} \tag{13}$$

$$= \frac{\mu(\mathcal{X}_\mathbf{A}) \cdot [\overline{\mu} \star \overline{\lambda}]}{\mu(\mathcal{X}_\mathbf{A}) \cdot [\overline{\mu} \star \overline{\lambda}](\mathcal{X}_{\mathbf{A} \cup \mathbf{B}})} \tag{14}$$

$$= \overline{\mu} \star \overline{\lambda}. \qquad \Box \tag{15}$$

Writing the base measures in the alternative notation, set $\mu = \nu G_1$ and $\lambda = \nu G_2$. Theorem 13 states that $\overline{\mu \star \lambda} = G_1 \star G_2$. Therefore, the Markov combination of $\nu G_1$ and $\nu G_2$ can be written $\nu(G_1 \star G_2)$.

### 4.3. Constructing the Hyper Dirichlet Process

We now apply the idea of Markov combinations to component Dirichlet processes. To do so, we simply form the Markov combination of the base measures.

**Theorem 14** *Let $G_1$ be a distribution on $(\mathcal{X}_\mathbf{A}, \mathcal{F}_\mathbf{A})$. Let $H_2$ be a distribution on $(\mathcal{X}_\mathbf{A}, \mathcal{F}_\mathbf{A})$ that is consistent with $H_1$. Set $H = H_1 \star H_2$. Let $Q \sim DP(\nu H_1)$, $R \sim DP(\nu H_2)$, and $\theta \sim DP(\nu H)$ be random probability measures. The following are true:*

*(i) $\theta_\mathbf{A} \stackrel{d}{=} Q$*

*(ii) $\theta_\mathbf{B} \stackrel{d}{=} R$*

PROOF: The proposition follows from Theorem 5. $\Box$

Note that the non-parametric approach is actually simpler than the parametric approach in one sense. The hyper inverse Wishart is a generalization of the inverse Wishart to incomplete graphs. However, the distribution of $\theta$ in Theorem 14 actually *is* a Dirichlet process with a Markov base measure. Therefore, previous results regarding the Dirichlet process also apply to the hyper Dirichlet process. Most importantly, we know that the prior law is a stick-breaking prior.

The hyper inverse Wishart is so-called because it is an example of a hyper Markov prior (Dawid and Lauritzen, 1993).

**Definition 15** HYPER MARKOV. *Consider an undirected graph, $\mathcal{G}$. Let $\theta \sim \mathcal{L}$ be a random probability measure over $\mathcal{X}$. We say that $\mathcal{L}$ is (weak) hyper Markov on $\mathcal{G}$ if $\mathcal{L}$ is concentrated on $\mathscr{M}(\mathcal{G})$, and $\theta_\mathbf{A} \perp\!\!\!\perp \theta_\mathbf{B} | \theta_\mathbf{C}$ whenever $\mathbf{C}$ separates $\mathbf{A}$ and $\mathbf{B}$.*

It is tempting to refer to the distribution of $\theta$ in Theorem 14 as a hyper Dirichlet process; however, the given conditions are not sufficient to ensure that the process is hyper Markov. The next task is to discover the appropriate conditions under which the hyper Markov property holds. Let $\mathcal{L} = DP(\nu H)$ be a Dirichlet process law. Let $\mathcal{G}$ be any graph consisting of two cliques, $\mathbf{A}$ and $\mathbf{B}$, with separator $\mathbf{C}$. Using the stick-breaking construction, let $\vec{w} = (w_i : i \in \mathbb{N})$ be the random weights and $\vec{Z} = (Z_i : i \in \mathbb{N})$ be the atoms, which are iid observations from $H$. We use $Z_{i\mathbf{\Gamma}}$ to denote the components of $Z_i$ belonging to a set $\mathbf{\Gamma}$. For example, the marginal of $\theta$ over $\mathbf{A}$ is $\theta_\mathbf{A} = \sum_{i \in N} w_i \delta_{Z_{i\mathbf{A}}}$.

Obviously, one condition for hyper Markovity is that $H$ is a Markov measure. If $H$ is not Markov, then for $Z_i \sim H$, it is not true that $Z_{i\mathbf{A}} \perp\!\!\!\perp Z_{i\mathbf{B}} | Z_{i\mathbf{C}}$. As a result, $\theta_\mathbf{A} \not\perp\!\!\!\perp \theta_\mathbf{B} | \theta_\mathbf{C}$.

$H \in \mathscr{M}(\mathcal{G})$ is a necessary condition, but it is not sufficient. This is because knowledge of $\theta_\mathbf{B}$ contains information about the distribution of weights at each atom. We must ensure that $\theta_\mathbf{C}$ contains the information as well. To see this, consider an example for which $H_\mathbf{C}$ is a point mass. For $H \in \mathscr{M}(\mathcal{G})$, this implies that $H_\mathbf{A} \perp\!\!\!\perp H_\mathbf{B}$. Further suppose that $H_\mathbf{A}$ and $H_\mathbf{B}$ are *not* point masses. In this case, $\theta_\mathbf{B}$ implies certain constraints on $\vec{w}$. For example, if each $Z_{i\mathbf{B}}$ is distinct, then the mass at each atom determines the random weights modulo permutation. Therefore, the second condition for $\mathcal{L}$ to be hyper Markov is that $\theta_\mathbf{B}$ contains no information about $\vec{w}$ that is not contained by $\theta_\mathbf{C}$. To begin, we use the condition expressed in the next theorem. The condition is sufficient, but more restrictive than necessary.

**Theorem 16** *Let $H$ be a base measure on $\mathcal{X}_{\mathbf{A} \cup \mathbf{B}}$. Let $\mathbf{C} = \mathbf{A} \cap \mathbf{B}$. Set $\mathcal{L} = DP(\nu H)$ for some $\nu > 0$. Then $\mathcal{L}$ is hyper Markov on $\mathcal{X}_{\mathbf{A} \cup \mathbf{B}}$ if the following conditions hold:*

1. *$H$ is a Markov measure*
2. *Refinement Condition:*

$$Z_{i\mathbf{C}} = Z_{j\mathbf{C}} \Rightarrow Z_{i\mathbf{B}} = Z_{j\mathbf{B}} \ \text{a.s.}[H].$$

PROOF: Define $A' = \mathbf{A} \setminus \mathbf{C}$ and $B' = \mathbf{B} \setminus \mathbf{C}$. Note that $\mathbf{B} = \mathbf{C} \cup \mathbf{B}'$, so that $Z_{i\mathbf{B}} = Z_{j\mathbf{B}} \Rightarrow (Z_{i\mathbf{C}}, Z_{i\mathbf{B}'}) = (Z_{j\mathbf{C}}, Z_{i\mathbf{B}'})$. In other words, the refinement condition can be expressed equivalently as an "if and only if" statement:

$$Z_{i\mathbf{C}} = Z_{j\mathbf{C}} \iff Z_{i\mathbf{B}} = Z_{j\mathbf{B}} \text{ a.s.}[H].$$

Consider $\theta \sim \mathcal{L}$. The hyper Markov property has two conditions:

1. $\mathbb{P}(\theta \in \mathscr{M}(\mathcal{G})) = 1$, and
2. $\theta_{\mathbf{A}} \perp\!\!\!\perp \theta_{\mathbf{B}} | \theta_{\mathbf{C}}$.

The first condition follows from the refinement condition. Let $x = (x_{\mathbf{A}'}, x_{\mathbf{C}}, x_{\mathbf{B}'})$ be any point in $\mathcal{X}$ such that $\theta_{\mathbf{B}}(x) > 0$. That is, there exists some $i$ such that $Z_{i\mathbf{B}} = x_{\mathbf{B}}$. By the refinement condition, $Z_{j\mathbf{C}} = Z_{i\mathbf{C}} = x_{\mathbf{C}}$ if and only if $Z_{j\mathbf{B}} = Z_{i\mathbf{C}} = x_{\mathbf{B}}$. Hence, $\{j : Z_{j\mathbf{C}} = x_{\mathbf{C}}\} = \{j : Z_{j\mathbf{B}} = x_{\mathbf{B}}\}$. Using the stick-breaking representation, we write the distribution of $X_{\mathbf{A}}|X_{\mathbf{B}}$.

$$\theta_{\mathbf{A}|\mathbf{B}}(x_{\mathbf{A}}|X_{\mathbf{B}} = x_{\mathbf{B}}) = \frac{\sum_{i: Z_{i\mathbf{B}} = x_{\mathbf{B}}} w_i 1_{\{Z_{i\mathbf{A}} = x_{\mathbf{A}}\}}}{\sum_{i: Z_{i\mathbf{B}} = x_{\mathbf{B}}} w_i} \tag{16}$$

$$= \frac{\sum_{i: Z_{i\mathbf{C}} = x_{\mathbf{C}}} w_i 1_{\{Z_{i\mathbf{A}} = x_{\mathbf{A}}\}}}{\sum_{i: Z_{i\mathbf{C}} = x_{\mathbf{C}}} w_i} \tag{17}$$

$$= \theta_{\mathbf{A}|\mathbf{C}}(x_{\mathbf{A}}|X_{\mathbf{C}} = x_{\mathbf{c}}). \tag{18}$$

Therefore, $\theta \in \mathscr{M}(\mathcal{G})$.

It remains to show that $\theta_{\mathbf{A}} \perp\!\!\!\perp \theta_{\mathbf{B}} | \theta_{\mathbf{C}}$. We begin by writing the marginals of $\theta$ using the stick-breaking representation. Let $\Gamma$ be any subset of $\mathbf{V}$.

$$\theta_{\mathbf{\Gamma}} = \sum_{i \in \mathbb{N}} w_i \delta_{Z_{i\mathbf{\Gamma}}}. \tag{19}$$

Let $\vec{Z}_{\mathbf{\Gamma}}^* = \{Z_{i\mathbf{\Gamma}}\}$ be the set of unique occurrences among the random atoms. We refer to an element of this set using an arbitrary index, $Z_{i\mathbf{\Gamma}}^*$. Let $m_{i\mathbf{\Gamma}}^*$ be the total mass at that atom;

$$m_{i\mathbf{\Gamma}}^* = \sum_{j: Z_{j\mathbf{\Gamma}} = Z_{i\mathbf{\Gamma}}^*} w_j = \theta_{\mathbf{\Gamma}}(Z_{i\mathbf{\Gamma}}^*). \tag{20}$$

Note that $\vec{Z}_{\mathbf{\Gamma}}^*$ is the support of the $\theta_{\mathbf{\Gamma}}$, and $\vec{m}_{\mathbf{\Gamma}}^*$ is the mass at each point in the support. Thus, there is a bijection between the measure $\theta_{\mathbf{\Gamma}}$ and the set $\{\vec{Z}_{\mathbf{\Gamma}}^*, \vec{m}_{\mathbf{\Gamma}}^*\}$. That is to say, both are completely identified if at least one is known. The immediate result is that conditioning on one (or both) is equivalent to conditioning on the other (or one of them).

Continue by partitioning the support into two sets. Define $\vec{Z}_{\mathbf{\Gamma}}^+ = \{Z_{i\mathbf{\Gamma}}^* : H_{\mathbf{\Gamma}}(Z_{i\mathbf{\Gamma}}^*) > 0\}$ and $\vec{Z}_{\mathbf{\Gamma}}^0 = \{Z_{i\mathbf{\Gamma}}^* : H_{\mathbf{\Gamma}}(Z_{i\mathbf{\Gamma}}^*) = 0\} = \vec{Z}_{\mathbf{\Gamma}}^* \setminus \vec{Z}_{\mathbf{\Gamma}}^+$. In other words $\vec{Z}_{\mathbf{\Gamma}}^+$ is the set of support points with strictly positive mass and $\vec{Z}_\Gamma^0$ is the set of points that are in the support but have probability zero. Again, we specify a particular element in either set with an arbitrary index, e.g. $Z_{i\Gamma}^+$. Partition $\vec{m}_\Gamma^*$ in the same way. This yields,

$$\vec{m}_\Gamma^+ = \{m_{i\Gamma}^* : H_\Gamma(Z_{i\Gamma}^*) > 0\} = \{m_{i\Gamma}^* : Z_{i\Gamma}^* \in \vec{Z}_\Gamma^+\} = \{\theta_\Gamma(Z_{i\Gamma}^+)\}. \qquad (21)$$

We stipulate that the index is consistent with $\vec{Z}_\Gamma^+$ so that $m_{i\Gamma}^+ = \theta_\Gamma(Z_{i\Gamma}^+)$. Denote the other set in this partition by $\vec{m}_\Gamma^0 = \vec{m}_\Gamma^* \setminus \vec{m}_\Gamma^+$, where $m_{i\Gamma}^0 = \theta_\Gamma(Z_{i\Gamma}^0)$. Separate the sum in Equation 19 using this partition.

$$\theta_\Gamma = \sum_{i=1}^{N_\Gamma^+} m_{i\Gamma}^+ \delta_{Z_{i\Gamma}^+} + \sum_{i=1}^{N_\Gamma^0} m_{i\Gamma}^0 \delta_{Z_{i\Gamma}^0} , \qquad (22)$$

where $N_\Gamma^\cdot = |\vec{Z}_\Gamma^\cdot|$. Note that $Z_\Gamma^+$, has a degenerate distribution. If $H_\Gamma(\gamma) > 0$, then with probability 1, $\gamma$ will occur infinitely often in $\vec{Z}_\Gamma$. Therefore, $\vec{Z}_\Gamma^+ = \{z_\Gamma : H_\Gamma(z_\Gamma) > 0\}$ almost surely. Since $H$ is known, the sets of summation in Equation 22 are fully identified by $\{\vec{Z}_\Gamma^*, \vec{m}_\Gamma^*\}$. It follows that conditioning on $\theta_\Gamma$ is equivalent to conditioning on the quartet $\{\vec{Z}_\Gamma^+, \vec{m}_\Gamma^+, \vec{Z}_\Gamma^0, \vec{m}_\Gamma^0\}$. We will now show that under the refinement condition, $\vec{Z}_B^+, \vec{m}_B^+$, and $\vec{m}_B^0$ are fully identified from $\theta_C$. With that fact, showing $\theta_A \perp\!\!\!\perp \theta_B | \theta_C$ is equivalent to showing $\theta_A \perp\!\!\!\perp \vec{Z}_B^0 | \theta_C$.

By the refinement condition,

$$m_{i\mathbf{C}}^+ = \sum_{j:Z_{j\mathbf{C}}=Z_{i\mathbf{C}}^+} w_j = \sum_{j:Z_{j\mathbf{B}}=Z_{i\mathbf{B}}^+} w_j = m_{i\mathbf{C}}^+. \qquad (23)$$

A similar equation shows $m_{i\mathbf{C}}^0 = m_{i\mathbf{B}}^0$. Therefore, $(\vec{m}_\mathbf{C}^+, \vec{m}_\mathbf{C}^0) = (\vec{m}_\mathbf{B}^+, \vec{m}_\mathbf{B}^0)$.

We now show that $\vec{Z}_\mathbf{B}^+$ is a function of $\vec{Z}_\mathbf{C}^+$. This fact ensures that $\vec{Z}_\mathbf{B}^+$ is fully identified by $\vec{Z}_\mathbf{C}^+$ and therefore is conditionally independent of anything given $\vec{Z}_\mathbf{C}^+$. One consequence of the refinement condition is that if $H_\mathbf{C}(c) > 0$, then there exists $B(c)$ such that $H_{\mathbf{B}|\mathbf{C}}(B(c)|c) = 1$. This follows from a simple proof by contradiction. If $H_{\mathbf{B}|\mathbf{C}}(\cdot|c)$ is not a point distribution, than either every point has probability 0 (as in a continuous distribution), or there is some point with positive probability strictly less than 1. We will see that neither of these can be true and conclude the conditional is indeed a point distribution.

Suppose $H_{\mathbf{B}|\mathbf{C}}(\cdot|c)$ has measure zero everywhere. With probability $H_\mathbf{C}(c)^2 > 0$, the event $Z_{1\mathbf{C}} = Z_{2\mathbf{C}}$ will occur. However, $Z_{1\mathbf{B}} \neq Z_{2\mathbf{B}}$ almost surely. Therefore, the refinement condition fails with probability at least $H_\mathbf{C}(c)^2 > 0$. Now suppose there exists $b$ such that $0 < H_{\mathbf{B}|\mathbf{C}}(b|c) < 1$. Then with probability $H_\mathbf{C}(c)^2 H_{\mathbf{B}|\mathbf{C}}(b|c)(1 - H_{\mathbf{B}|\mathbf{C}}(b|c)) > 0$, the events $Z_{1\mathbf{C}} = c = Z_{2\mathbf{C}}$ and $Z_{1\mathbf{B}} = b \neq Z_{2\mathbf{B}}$ will occur. Thus, the refinement condition fails with positive probability. By these two contradictions, we see that $H_{\mathbf{B}|\mathbf{C}}(\cdot|c)$ must be a point distribution if $H_\mathbf{C}(c) > 0$. We denote the point of concentration by $B(c)$. Clearly, $c \in \vec{Z}_\mathbf{C}^+$ implies that $B(c) \in \vec{Z}_\mathbf{B}^+$. Furthermore, every element of $\vec{Z}_\mathbf{B}^+ = B(c)$ for some $c \in \vec{Z}_\mathbf{C}^+$. This follows from the fact that $\mathbf{C} \subseteq \mathbf{B}$, so $H_\mathbf{C}(Z_{i\mathbf{C}}) = 0 \Rightarrow H_\mathbf{B}(Z_{i\mathbf{C}}, Z_{i\mathbf{B}'}) = 0$. Therefore, $\vec{Z}_\mathbf{B}^+ = g(\vec{Z}_\mathbf{C}^+) = \{B(c) : c \in \vec{Z}_\mathbf{C}^+\}$ almost surely.

We have shown that conditioning on $\theta_.$ is equivalent to conditioning on $\{\vec{Z}_.^+, \vec{m}_.^+, \vec{Z}_.^0, \vec{m}_.^0\}$. Furthermore, we have that $(\vec{Z}_\mathbf{B}^+, \vec{m}_\mathbf{B}^+, \vec{m}_\mathbf{B}^0) = (g(\vec{Z}_\mathbf{C}^+), \vec{m}_\mathbf{C}^+, \vec{m}_\mathbf{C}^0)$. This provides an equivalent condition for the independence property that we want to show. That is, $\theta_\mathbf{A} \perp\!\!\!\perp \theta_\mathbf{B} | \theta_\mathbf{C}$ if and only if $\theta_\mathbf{A} \perp\!\!\!\perp \vec{Z}_\mathbf{B}^0 | \{\vec{Z}_\mathbf{C}^+, \vec{m}_\mathbf{C}^+, \vec{Z}_\mathbf{C}^0, \vec{m}_\mathbf{C}^0\}$. The remainder of this proof will show that the second property holds under the conditions of the theorem.

Begin by partitioning the atoms and weights as follows. Let $\hat{Z} = \{Z_i : Z_{i\mathbf{C}} \in \vec{Z}_\mathbf{C}^+\}$, and $\tilde{Z} = \{Z_i : Z_{i\mathbf{C}} \in \vec{Z}_\mathbf{C}^0\}$. Let $\hat{w} = \{w_i : Z_{i\mathbf{C}} \in \vec{Z}_\mathbf{C}^+\}$, and $\tilde{w} = \{w_i : Z_{i\mathbf{C}} \in \vec{Z}_\mathbf{C}^0\}$. As usual, for $\mathbf{\Gamma} \subseteq \mathbf{V}$, let $\hat{Z}_\mathbf{\Gamma}$ and $\tilde{Z}_\mathbf{\Gamma}$ denote that the elements are the components in $\mathbf{\Gamma}$. This partition is similar to, but different than, the partition defined earlier. $(\hat{Z}_\mathbf{\Gamma}, \tilde{Z}_\mathbf{\Gamma})$ depends on $H_\mathbf{C}$, whereas $(Z_\mathbf{\Gamma}^+, Z_\mathbf{\Gamma}^0)$ depends on $H_\mathbf{\Gamma}$. The goal, as above is to rewrite $Z_\mathbf{A}$ by partitioning it in a way that preserves the conditional independence structure. This structure is preserved if the partioning function is non-random. In other words, the atoms must be partitioned based on on a known event. When conditioning on $\theta_\mathbf{B}$, $\theta_\mathbf{B}$ and $\theta_\mathbf{C}$ are known, but $\theta_\mathbf{A}$ is unknown. Therefore, $(\hat{Z}_\mathbf{A}, \tilde{Z}_\mathbf{A})$ provides an observable partition of $\vec{Z}_\mathbf{A}$.

$$\theta_\mathbf{A} = \sum_{i=1}^{N_\mathbf{C}^+} \hat{w}_i \hat{Z}_{i\mathbf{A}} + \sum_{i=1}^{N_\mathbf{C}^0} \tilde{w}_i \tilde{Z}_{i\mathbf{A}}. \quad (24)$$

Note that $\tilde{Z}_\mathbf{C}$ is equivalent to $\vec{Z}_\mathbf{C}^0$ by definition, and $\tilde{Z}_\mathbf{B} = \vec{Z}_\mathbf{B}^0$ by the refinement condition. We proceed by showing that $\tilde{w}, \tilde{Z}_\mathbf{A}, \hat{w}$, and $\hat{Z}_\mathbf{A}$ are jointly independent of $\vec{Z}_\mathbf{B}^0$ given $\{\vec{Z}_\mathbf{C}^+, \vec{m}_\mathbf{C}^+, \vec{Z}_\mathbf{C}^0, \vec{m}_\mathbf{C}^0\}$. We can express $\vec{m}_\mathbf{C}^+$ as a function of $\hat{w}, \hat{Z}$, and $\vec{Z}_\mathbf{C}^+$, where

$$m_{i\mathbf{C}}^+ \stackrel{a.s.}{=} \sum_{j:\hat{Z}_j = Z_{i\mathbf{C}}^+} \hat{w}_j. \quad (25)$$

Furthermore, we have noted that $\vec{m}_\mathbf{C}^0 = \tilde{w}$. By the stick-breaking construction, $\tilde{Z} \perp\!\!\!\perp (\tilde{w}, \hat{w}, \hat{Z})$. Since $\vec{Z}_\mathbf{C}^+$ is known almost surely, it can be included in the independence property.

$$\tilde{Z} \perp\!\!\!\perp (\vec{Z}_\mathbf{C}^+, \hat{Z}, \hat{w}, \tilde{w}). \quad (26)$$

$\tilde{Z}$ is also independent of any function of the RHS of Equation 26. In particular,

$$\tilde{Z} \perp\!\!\!\perp (\vec{m}_\mathbf{C}^+, \vec{m}_\mathbf{C}^0, \vec{Z}_\mathbf{C}^+, \hat{Z}_\mathbf{A}, \hat{w}, \tilde{w}). \quad (27)$$

Repeating this argument on the LHS of Equation 26, we conclude

$$(\vec{Z}_\mathbf{B}^0, \tilde{Z}_\mathbf{A}, \vec{Z}_\mathbf{C}^0) \perp\!\!\!\perp (\vec{m}_\mathbf{C}^+, \vec{m}_\mathbf{C}^0, \vec{Z}_\mathbf{C}^+, \hat{Z}_\mathbf{A}, \hat{w}, \tilde{w}). \quad (28)$$

Since $\tilde{Z}_i \sim H \in \mathscr{M}(\mathcal{G})$, we can write $\vec{Z}_\mathbf{B}^0 \perp\!\!\!\perp \tilde{Z}_\mathbf{A} | \vec{Z}_\mathbf{C}^0$. Since all three of these are jointly independent of $(\hat{Z}_\mathbf{A}, \hat{w})$ and $(\vec{m}_\mathbf{C}^+, \vec{m}_\mathbf{C}^0, \vec{Z}_\mathbf{C}^+))$, it follows that

$$\vec{Z}_{\mathbf{B}}^0 \perp\!\!\!\perp (\hat{Z}_{\mathbf{A}}, \hat{w}, \tilde{Z}_{\mathbf{A}}, \tilde{w}) | (\vec{m}_{\mathbf{C}}^+, \vec{m}_{\mathbf{C}}^0, \vec{Z}_{\mathbf{C}}^+, \vec{Z}_{\mathbf{C}}^0). \tag{29}$$

Recall from Equation 24 that $\theta_{\mathbf{A}}$ is a function of $(\hat{Z}_{\mathbf{A}}, \hat{w}, \tilde{Z}_{\mathbf{A}}, \tilde{w})$. A discrete distribution can be defined by a set of points and the probability mass at each point. Therefore, $\theta_{\mathbf{A}}$ is a function of $(\hat{Z}_{\mathbf{A}}, \hat{w}, \tilde{Z}_{\mathbf{A}}, \tilde{w})$. It follows that,

$$\vec{Z}_{\mathbf{B}}^0 \perp\!\!\!\perp \theta_{\mathbf{A}} | (\vec{m}_{\mathbf{C}}^+, \vec{m}_{\mathbf{C}}^0, \vec{Z}_{\mathbf{C}}^+, \vec{Z}_{\mathbf{C}}^0). \tag{30}$$

Hence, by the above argument, it follows that $\theta_{\mathbf{A}} \perp\!\!\!\perp \theta_{\mathbf{B}} | \theta_{\mathbf{C}}$. We conclude that $\mathcal{L}$ is hyper Markov. □

Theorem 16 provides sufficient conditions for a Dirichlet process to be hyper Markov. Thus, under those conditions, we may safely call the Dirichlet process a *hyper* Dirichlet process. When $H$ satisfies the refinement condition, we will say that $\mathbf{B}$ is a refinement of $\mathbf{C}$ under sampling almost surely under measure $H$. It is a refinement in the following sense. Let $X_1, X_2, \ldots$ be an infinite iid sample from $H$. Form a partition of the natural numbers such that $i$ and $j$ are elements of the same set if and only if $X_{i\mathbf{C}} = X_{j\mathbf{C}}$. Call this partition $X(\mathbf{C})$. Define $X(\mathbf{B})$ by analogy. Under the refinement condition, $X(\mathbf{B})$ is almost surely a refinement of $X(\mathbf{C})$. We denote this relationship by $X(\mathbf{B}) \preceq X(\mathbf{C})$ $a.s.[H]$, omitting $H$ if the measure is contextually evident.

The refinement condition, as stated in Theorem 16 is sufficient, but it is stronger than necessary. By symmetry of conditional independence, $\theta_{\mathbf{B}} \perp\!\!\!\perp \theta_{\mathbf{A}} | \theta_{\mathbf{C}}$, even though no refinement condition is needed between $\mathbf{C}$ and $\mathbf{A}$. It may be necessary that at least one of the two refinements is present, but this has not been explored.

The hyper Dirichlet process defined on two cliques is an example of a hyper Markov combination, which is the analog of Markov combinations for prior laws. Consider two laws: $\mathcal{Q}$ for $\theta_{\mathbf{A}}$ and $\mathcal{R}$ for $\theta_{\mathbf{B}}$. We say that $\mathcal{Q}$ and $\mathcal{R}$ are *hyperconsistent* if the marginal laws for $\theta_{\mathbf{A} \cup \mathbf{B}}$ are equal. Under this condition, Dawid and Lauritzen (1993) show that there is a unique hyper Markov law $\mathcal{L}$ such that $\mathcal{L}_{\mathbf{A}} = \mathcal{Q}$, $\mathcal{L}_{\mathbf{B}} = \mathcal{R}$. This is called the hyper Markov combination and is denoted $\mathcal{L} = \mathcal{Q} \odot \mathcal{R}$.

As with Markov combinations, hyper Markov combinations are easily generalized to multiple cliques. Let $\mathcal{G}$ be a graph with perfectly ordered cliques $(\mathbf{C}_1, \ldots, \mathbf{C}_k)$. Suppose $\mathbf{C}_i$ is imbued with a prior law $\mathcal{G}_i$ and that the priors are all pairwise hyperconsistent. Let $\mathcal{L}_1 = \mathcal{G}_1$ and $\mathcal{L}_i = \mathcal{L}_{i-1} \odot \mathcal{G}_i$ for $i \geq 2$. Then $\mathcal{L} = \mathcal{L}_k$ is the unique hyper Markov prior satisfying $\mathcal{L}_{\mathbf{C}_i} = \mathcal{G}_i$. We call $\mathcal{L}$ the hyper Markov combination of $\mathcal{G}_1, \ldots, \mathcal{G}_k$. In general we may write $\odot(\mathcal{G}_1, \ldots, \mathcal{G}_k)$ with the understanding that the cliques are perfectly ordered and $\mathcal{G}_1, \ldots, \mathcal{G}_k$ are pairwise consistent.

The next definition generalizes the hyper Dirichlet Process to three or more cliques.

**Definition 17** HYPER DIRICHLET PROCESS. *Let $\mathcal{G}$ be a graph with a perfect ordering of cliques $\mathbf{C}_1, \ldots, \mathbf{C}_k$. Suppose that the $i^{\text{th}}$ clique has marginal distribution $G_i$ and that the marginals are pairwise consistent. Let $G = \star(G_1, \ldots, G_k)$. Further suppose that $\mathbf{C}_j$ or $\mathbf{H}_j$ is a refinement of $\mathbf{S}_j$ under sampling almost surely under $H$, where $\mathbf{H}_i$ is the $i^{\text{th}}$ history and $\mathbf{S}_i$ is the $i^{\text{th}}$ separator. Then*

$$\text{HDP}(\nu, G_1, \ldots, G_k) = \text{DP}(\nu G) \tag{31}$$

*is a hyper Dirichlet process prior.*

This hyper Dirichlet process defined in this way is guaranteed to be hyper Markov. Suppose $\mathcal{L} = \text{HDP}(\nu, G_1, \ldots, G_k)$. By Theorem 7, $\mathcal{L}_{\mathbf{C}_i} = \text{DP}(\nu G_i)$ for $i \geq 2$. Furthermore, it follows from the refinement conditions and Theorem 16 that $\mathcal{L}_{\mathbf{H}_{i-1} \cup \mathbf{C}_i} = \text{DP}(\nu G_{\mathbf{H}_{i-1}}) \odot \text{DP}(\nu G_i)$. Hence, Theorem 3.9 in Dawid and Lauritzen (1993) states that $\mathcal{L}$ is the almost-everywhere unique hyper Markov law such that $\mathcal{L}_{\mathbf{C}_i} = \text{DP}(\nu G_i)$.

The next theorem states that if the prior distribution of $\theta$ is a hyper Dirichlet process, then so is the posterior.

**Theorem 18** POSTERIOR HYPER DIRICHLET PROCESS. *Suppose $\mathcal{G}$ is a graph with a perfect ordering of cliques $\mathbf{C}_1, \ldots, \mathbf{C}_k$. Let $\mathcal{L} = \text{HDP}(\nu, H_1, \ldots, H_k)$ be a hyper Dirichlet process. Given $\theta \sim \mathcal{L}$, let $X_1, \ldots, X_n$ be an iid sample from $\theta$. Denote the marginal value of $X_j$ over the $i^{\text{th}}$ clique by $X_{ji}$. Then,*

$$\theta | X_1, \ldots, X_n \sim \text{HDP}(\nu, H'_1, \ldots, H'_k), \tag{32}$$

*where $H'_i = (\nu + n)^{-1}(\nu H + \sum_j \delta_{X_{ji}})$*

PROOF: By definition $\theta \sim \text{DP}(\nu H)$, where $H$ is the Markov combination of $H_1, \ldots, H_k$. Therefore, Theorem 7 states that $\theta | X_1, \ldots, X_n \sim \text{DP}(\nu' H')$, where $\nu' = \nu + n$ and $H' = (\nu + n)^{-1}(\nu H + \sum_j \delta_{X_j})$. Furthermore, the marginal distribution of $H'$ over the $i^{\text{th}}$ clique is $H'_i$. It follows that $H'$ is the unique Markov Combination of $H'_1, \ldots, H'_k$. It remains to show that $\text{DP}(\nu' H')$ is a hyper Markov law via the refinement condition.

For $1 < i \leq k$, let $\mathbf{S}_i$ be the $i^{\text{th}}$ separator and let $\mathbf{G}_i$ be the $i^{\text{th}}$ history. By Definition 17, $X(\mathbf{A}_i) \preceq X(\mathbf{S}_i)$ $a.s.[H]$, where $\mathbf{A}_i$ is either the $i^{\text{th}}$ clique or the $i^{\text{th}}$ history. Let $Z$ be such that $H'_\mathbf{S}(Z) > 0$. We will show that there exists some $a$, such that $H_{\mathbf{A}|\mathbf{S}}(\cdot|Z) = \delta_a(\cdot)$. From this it follows that $X(\mathbf{A}_i) \preceq X(\mathbf{S}_i)$ $a.s.[H']$. We consider the cases $H_\mathbf{S}(Z) > 0$ and $H_\mathbf{S}(Z) = 0$ separately.

Suppose $H_\mathbf{S}(Z) > 0$. By the refinement condition, there exists $a$ such that $H_{\mathbf{A}|\mathbf{S}}(\cdot|Z) = \delta_a(\cdot)$. Therefore, if $X_{j\mathbf{S}} = Z$, then $X_{j\mathbf{A}} = a$ with probability one. Hence,

$$
\begin{aligned}
H'_{\mathbf{A}|\mathbf{S}}(\cdot|Z) &= \frac{\nu H_\mathbf{S}(Z) H_{\mathbf{A}|\mathbf{S}}(\cdot|Z) + \sum_i \delta_{X_{j\mathbf{A}}}(Z, \cdot)}{\nu H_\mathbf{S}(Z) + \sum_j \delta_{X_{j\mathbf{S}}}(Z)} \tag{33} \\
&\stackrel{a.s.}{=} \frac{\nu H_\mathbf{S}(Z) \delta_a(\cdot) + \sum_j \delta_a(\cdot)}{\nu H_\mathbf{S}(Z) + \sum_j \delta_a(\cdot)} \tag{34}
\end{aligned}
$$

$$= \delta_a(\cdot). \tag{35}$$

On the other hand, suppose $H_{\mathbf{S}}(Z) = 0$. Since $H'_{\mathbf{S}}(Z) > 0$, there is some $i$ such that $X_{i\mathbf{S}} = Z$. Furthermore, with probability one it holds that $X_{j\mathbf{S}} \neq Z$ for all $j \neq i$. Therefore,

$$H'_{\mathbf{A}|\mathbf{S}}(\cdot|Z) = \frac{\nu H_{\mathbf{S}}(Z) H_{\mathbf{A}|\mathbf{S}}(\cdot|Z) + \sum_i \delta_{X_{j\mathbf{A}}}(Z, \cdot)}{\nu H_{\mathbf{S}}(Z) + \sum_j \delta_{X_{j\mathbf{S}}}(Z)} \tag{36}$$

$$\stackrel{a.s.}{=} \frac{\delta_{X_{j\mathbf{A}}}(\cdot)}{\delta_{X_j\mathbf{S}}(Z)} \tag{37}$$

$$= \delta_{X_{j\mathbf{A}}}(\cdot). \tag{38}$$

The last equation holds because $\mathbf{A} \subset \mathbf{S}$, so $\delta_{X_{j\mathbf{A}}}(\cdot) = 1$ is a stronger condition than $\delta_{Xj\mathbf{S}}(Z) = 1$.

From these two cases, we see that if $H'_{\mathbf{S}}(Z) > 0$ then $H'_{\mathbf{A}|\mathbf{S}}(\cdot|Z) = \delta_a(\cdot)$ for some $a$ which depends on $Z$. Ergo, for $X, Y \sim H'$, $X_{\mathbf{S}} = Y_{\mathbf{S}}$ implies $Z_{\mathbf{A}} = Y_{\mathbf{A}}$ almost surely and $DP(\nu H)$ is a hyper Markov measure. □

## 5. Applications

At first glance, the refinement condition seems unduly restrictive. However, it allows one to use hyper Dirichlet processes in most areas that have benefited from Dirichlet processes. For the many applications that use a continuous base measure, the random atoms are distinct with probability one, so the refinement condition is trivial. Furthermore, Theorem 18 states that the posterior will also be hyper Markov. Thus, the hyper Dirichlet process can replace the Dirichlet process in applications requiring a posterior estimate, as in MCMC.

The hyper Dirichlet process provides a non-parametric alternative to the hyper Markov laws currently used for problems such as covariance selection for fitting graphical Gaussian models. Most previous work has focused on hyper inverse Wishart priors over $Q_\mathcal{G}$ (Roverato and Whittaker, 1998; Giudici and Green, 1999; Letac and Massam, 2007). The hyper inverse Wishart is conjugate to the hyper Wishart, which is the distribution of the restricted maximum likelihood estimate for the Gaussian problem with known mean. However, covariance selection could also be achieved using hyper Dirichlet processes. This provides a non-parametric prior for the covariance matrix. This could be advantageous if we wanted to remove the Gaussian assumption. A logical choice is to specify a hyper inverse Wishart for the base measure, setting the parameters by empirical Bayes. The precision parameter controls the concentration of the prior around this measure. Thus, it can loosely be considered a measure of confidence in the Gaussian assumption.

The real power of Dirichlet processes is in modeling mixture distributions. Suppose $X_1, \ldots, X_n$ are observations from some family parameterized by $\pi$. If we allow $\pi_i$ to be different for each observation, then the result is a mixture of distributions. The number of parameters increases with $n$, which necessitates placing some prior $G$ on the distribution of $\pi$ in order to fit the mixture. If the prior is unknown, it can be modeled with a Dirichlet process. For example, Escobar and West (1995) develop a Gibbs sampler to estimate the distribution of parameters of a Gaussian mixture model. In general, a Dirichlet mixture is a hierarchical model expressed as

$$\theta \sim DP_\alpha$$
$$\pi_1, \ldots, \pi_n | \theta \sim \theta$$
$$X_i | \pi_1, \ldots, \pi_n \sim f(X|\pi_i)$$

From Theorem 7, the posterior for $\theta$ given $\pi_1, \ldots, \pi_{n-1}$ is a Dirichlet process with base measure $\alpha(n-1) = a/(a+n-1)\overline{\alpha} + \sum_{i=1}^{n-1} 1/(a+n-1)\delta_{\pi_i}$, where $a = \alpha(\mathcal{X})$. Furthermore, this is also the distribution of $\pi_n | \pi_1, \ldots, \pi_{n-1}$. Thus, with positive probability, $\pi_n$ will be a previous value of $\pi_i$. Otherwise, it is drawn from $\overline{\alpha}$. As a result there will be $k \leq n$ unique values of $\pi_i$. This induces a latent class model for $X$ in which each class is defined by a shared value of $\pi_i$. The observations are conditionally independent given this latent class. A key feature of the Dirichlet process is that the number of latent classes is estimated. It is clear from the form of $\alpha(n)$ that this estimate is influenced by $a = \alpha(\mathcal{X})$. When $a$ is large, new values of $pi_i$ will often be drawn from $\alpha$. Contrarily, when $a$ is small, $\pi_i$ will more often be drawn from the previous values.

This is a natural setting for hyper Dirichlet processes. If $X$ is a multivariate random variable, then $\theta$ can be a hyper Dirichlet process for some graph. Once again the observations will be conditionally independent given their latent class. Furthermore, the components of $X$ will have the independence structure specified by the graph. For example, Escobar and West (1995) develop an MCMC algorithm for estimating Dirichlet mixture of Gaussian distributions. This can be extended to a mixture of the family $\mathcal{N}_\mathcal{G}$ of Gaussians that are Markov on $\mathcal{G}$ by restricting the base measure to $Q_\mathcal{G}$. For each update, the posterior of $\theta$ will be a hyper Markov prior.

As a final note, we point out that some of the results apply to other stick-breaking measures. Notably, Theorem 16 did not rely on the distribution of the random weights. Therefore, the same conditions imply that any stick-breaking measure is hyper Markov. That is, if the $H$ is Markov and the refinement condition holds, then a stick-breaking prior whose atoms have distribution $H$ is a hyper Markov law. Whether or not the posterior is also hyper Markov depends on how the measure is updated. For the Dirichlet process, the posterior update mechanism ensures a hyper Markov posterior as long as the prior is hyper Markov.

## 6. Acknowledgements

The author would like to thank Profs. S. Fienberg, A. Rinaldo, C. Schafer, and C. Shalizi for their guidance in developing this paper.

## Appendix A: Appendix: Working with Non-consistent Base Measures

We constructed the hyper Dirichlet process by combining base measures on each clique into a base measure on the entire graph. As we have seen, the end result is simply another Dirichlet process. In other words, the hyper Dirichlet process is simply a Dirichlet process that is hyper Markov on the graph of interest. The benefit is that the elicitation of the base measure is simplified by only considering a subset of variables at a time. In the current paper, we considered the case in which the individual base measures are consistent (i.e. they agree about the intersection.) However, we may not be able to guarantee consistency, especially if component base measures are elicited from different experts or models. How can this be resolved?

Consider a two cliques, $\mathbf{A}$ and $\mathbf{B}$ with intersection $\mathbf{C}$. Let $\mu$ be a measure on $\mathcal{X}_\mathbf{A}$ and $\lambda$ be a measure on $\mathcal{X}_\mathbf{B}$. The two measures are consistent if the marginals over $\mathcal{X}_\mathbf{C}$ are equal. In Section 4.2 we stated that this can be expressed as two simultaneous conditions: (1) the marginals must be proportional, and (2) the marginals must have the same scale.

Suppose only the first condition holds. Recall that $\mu(\mathcal{X}_\mathbf{A})$ represents the prior sample size. The interpretation is that there is more prior information about one clique than there is about the other clique. Equation 11 shows that there is no Markov Combination of $\mu$ and $\lambda$. That is, $\mu \star \lambda$ does not exist. Thus, there is no hyper Dirichlet process with those marginals. Fortunately, it is still possible to generate an "almost appropriate" random distribution. This is possible because we only need $\overline{\mu \star \lambda}$. By the commutative property, we can use $\overline{\mu} \star \overline{\lambda}$ instead. This is well-defined since the first property ensures that $\overline{\mu}$ and $\overline{\lambda}$ are consistent. On the other hand, the difference in scale must be resolved if a problem requires having a well-defined prior or posterior. The simplest way to achieve this is to scale one measure up or down to match the other. Additionally, any convex combination of $\mu\mathcal{X}_\mathbf{A}$ and $\lambda\mathcal{X}_\mathbf{B}$ could be a logical choice. The most conservative choice would be to set $\mu(\mathcal{X}_\mathbf{A}) = \lambda(\mathcal{X}_\mathbf{B}) = \min\{\mu(\mathcal{X}_\mathbf{A}), \lambda(\mathcal{X}_\mathbf{B})\}$. If resources are sufficient, both scales could be used and the results compared. This would reveal how sensitive the outcome is to the scale of the base measures.

Now suppose that the first condition does not hold. The interpretation is that we have conflicting prior information. Once again, $\mu \star \lambda$ does not exist. Furthermore, $\overline{\mu} \star \overline{\lambda}$ does not exist either, so it is not possible to use the same method to generate random distributions. In order to find a base measure, one or both distributions must be changed. There are several natural ways to do this. Let $U \subseteq \mathcal{X}_\mathbf{A}, V \subseteq \mathcal{X}_\mathbf{B}, W \subseteq \mathcal{X}_\mathbf{C}$.

1. Choose one base measure, and complete the distribution via conditioning.

$$\alpha^\mathbf{A}(U \times V \times W) = \mu(U \times W)\lambda(V|W). \tag{39}$$
$$\alpha^\mathbf{B}(U \times V \times W) = \mu(U|W)\lambda(V \times W). \tag{40}$$

2. Calculate a weighted average.

$$\alpha^w(U \times V \times W) = \gamma \alpha^{\mathbf{A}} + (1-\gamma)\alpha^{\mathbf{B}}, \quad \gamma \in [0,1]. \tag{41}$$

If there is no reason to choose one prior over the other, $\gamma = 1/2$ is appropriate. An interesting choice of $\gamma$ is $\mu(\mathcal{X}_{\mathbf{A}})/(\mu(\mathcal{X}_{\mathbf{A}}) + \lambda(\mathcal{X}_{\mathbf{B}}))$. This gives more weight toward the prior with more information.

3. Minimize the summed KL-divergence. Let $\mu_{\mathbf{C}}$ and $\lambda_{\mathbf{C}}$ be the marginals over $\mathcal{X}_{\mathbf{C}}$.

$$\alpha_{\mathbf{C}} = \underset{\nu}{\operatorname{argmin}} \int_{\mathcal{X}_{\mathbf{C}}} \overline{\mu}_{\mathbf{C}} \frac{\mu_{\mathbf{C}}}{\nu} + \int_{\mathcal{X}_{\mathbf{C}}} \overline{\lambda}_{\mathbf{C}} \frac{\lambda_{\mathbf{C}}}{\nu}. \tag{42}$$

$$\alpha(U \times V \times W) = \mu(U|W)\alpha_{\mathbf{C}}(W)\lambda(V|W). \tag{43}$$

More work is needed to test these candidate solutions to form good recommendations about their use.